\documentclass[11pt,a4,reqno]{amsart}
\usepackage{amssymb,latexsym}
\usepackage{graphicx}
\usepackage{amssymb}
\usepackage{url}		

\theoremstyle{plain}
\numberwithin{equation}{section}
\newtheorem{thm}{Theorem}[section]

\newtheorem{proposition}[thm]{Proposition}

\newtheorem{conjecture}{Conjecture}

\usepackage{mathrsfs}
\usepackage{diagbox}

\allowdisplaybreaks

\DeclareMathOperator{\sgn}{sgn}

\usepackage[font=small]{caption}

\usepackage{amsaddr}

\begin{document}
\setcounter{page}{1}

\title{Pi Visits Manhattan}

\author{Michelle Rudolph-Lilith}

\address{Unit\'e de Neurosciences, Information et Complexit\'e, CNRS\\
1 Ave de la Terrasse, 91198 Gif-sur-Yvette, France}

\email{rudolph@unic.cnrs-gif.fr}


\begin{abstract}
Is it possible to draw a circle in Manhattan, using only its discrete network of streets and boulevards? In this study, we will explore the construction and properties of circular paths on an integer lattice, a discrete space where the distance between two points is not governed by the familiar Euclidean metric, but the Manhattan or taxicab distance, a metric linear in its coordinates. In order to achieve consistency with the continuous ideal, we need to abandon Euclid's very original definition of the circle in favour of a parametric construction. Somewhat unexpectedly, we find that the Euclidean circle's defining constant $\pi$ can be recovered in such a discrete setting.
\end{abstract}



\maketitle


\section{Taxicabs in Manhattan.}

\begin{figure}[h!]
\centering
\includegraphics[width=\textwidth]{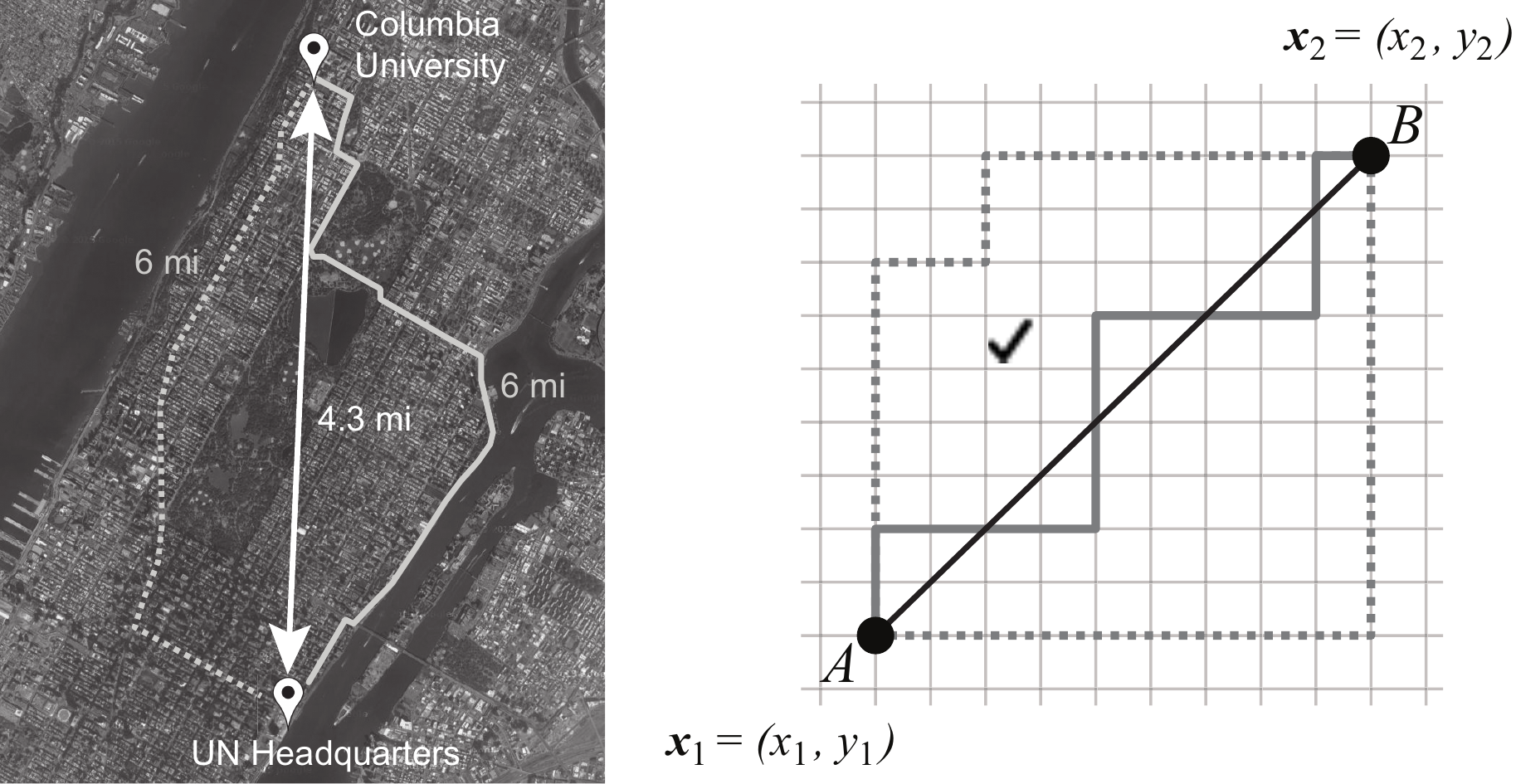}
\caption{\label{Fig_1}
Manhattan provides a prime example for the distance between two points in a discrete world. While the familiar Euclidean distance between two locations, here Columbia University and the UN Headquarters, is a mere 4.3 miles, the true distance a traveller will have traversed across Manhattan's busy streets is about 6 miles (left). Even more surprising, while on the Euclidean plane there is only one shortest path between two locations $A$ and $B$ (right; black), there is a multitude of shortest pathways when constrained to a discrete grid of streets and boulevards (gray).}
\end{figure}

A naive look at the map of Manhattan suggests that the distance between, for instance, Columbia University, located at 116th St \& Broadway, and the Headquarters of the United Nations at 1st Ave \& E 46th St is a mere 4.3 miles. But after leaving a taxicab, the attentive passenger will have spent about 6 miles traveling through the buzzing metropole (Fig.~\ref{Fig_1}). Slightly puzzled and with an emerging feeling of being cheated on, the passenger's mounting outrage will be met by a smiling cabdriver's shrug and cheerful reminder that ``this is New York''. 

Upon closer inspection, it will quickly become clear that, irrespective of whether the taxicab chooses to take on a rush-houred  East River Drive, the historical Broadway or whatever well-kept secret shortcut, these 6 miles distance are hardly to beat. The reason for this is the arrangement of streets and boulevards in Manhattan itself, which resembles a rather regular grid. On such a grid, the Euclidean distance with which we are all so intimately familiar with, has no say. Instead, a special case of distance measure, or metric, introduced more than a hundred years ago by the Russian-born mathematician and teacher of Albert Einstein, Hermann Minkowski, rules over Manhattan \cite{Minkowski67}. But it was not until the early 50ies of the last century that this metric obtained its now commonly used name ``taxicab'' or ``Manhattan distance'' \cite{Menger52}. 

On a mathematical more rigorous level, the Manhattan distance $\lVert . \rVert_1$ is called $\ell^1$-norm, and defined as the sum of the absolute value of the difference between the coordinates of two points (Fig.~\ref{Fig_1}, right). Restricting to the 2-dimensional Euclidean plane, and denoting with $\boldsymbol{x}_1 = (x_1,y_1)$ and $\boldsymbol{x}_2 = (x_2,y_2)$ two points $\boldsymbol{x}_1, \boldsymbol{x}_2 \in \mathbb{R}^2$, we have
\begin{equation}
\label{Eq_l1Norm}
\lVert \boldsymbol{x}_1 - \boldsymbol{x}_2 \rVert_1 = |x_1 - y_1| + |x_2 - y_2|.
\end{equation}
In contrast, the familiar Euclidean distance $\lVert . \rVert_2$, or $\ell^2$-norm, of these two points is defined as
\begin{equation}
\label{Eq_l2Norm}
\lVert \boldsymbol{x}_1 - \boldsymbol{x}_2 \rVert_2 = \sqrt{ (x_1 - y_1)^2 + (x_2 - y_2)^2}.
\end{equation}

Although (\ref{Eq_l1Norm}) looks somewhat simpler, as it is a measure of distance linear in its coordinates, it does not fall short of surprises. One of these unexpected properties we encountered already above in our trip through Manhattan. Whereas there is only one shortest path between two points on the Euclidean plane when measuring the distance with the $\ell^2$-norm, using the Manhattan distance does no longer provide a unique optimal path (Fig.~\ref{Fig_1}). In fact, in $\mathbb{R}^2$, we end up with a non-countable infinitude of shortest paths between any two points. Even restricting to a discrete lattice, a case resembling the situation in Manhattan, we still have multiple shortest pathways when their length is measured using the $\ell^1$-norm.

\begin{figure}[t!]
\centering
\includegraphics[width=\textwidth]{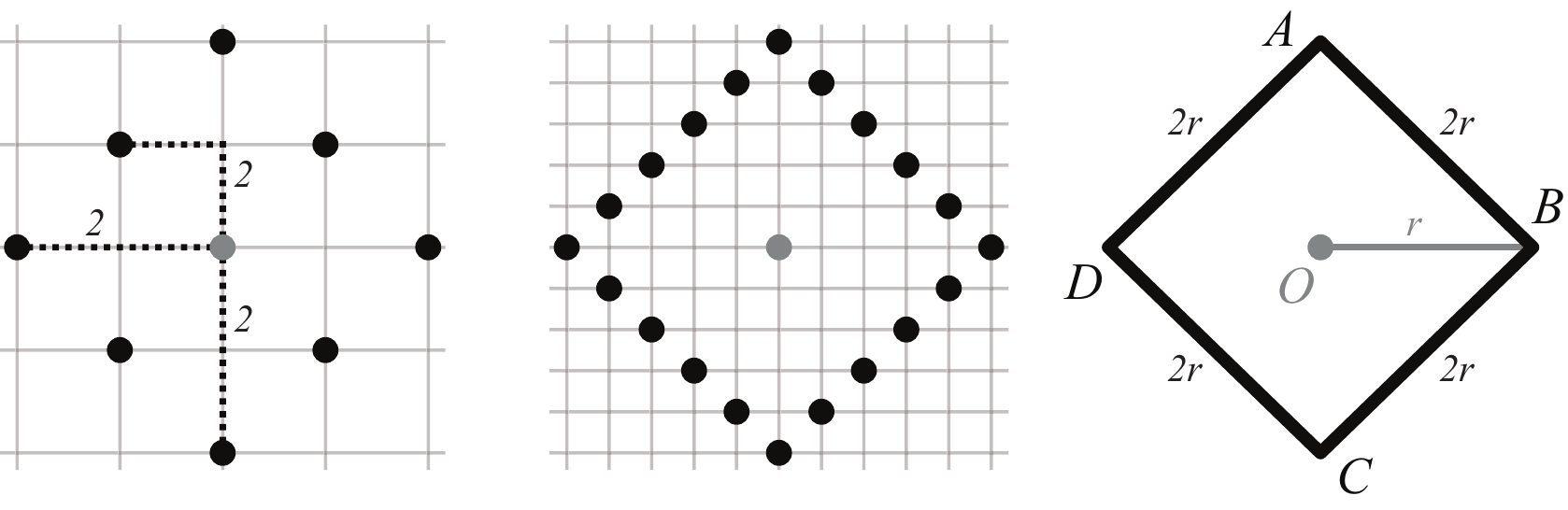}
\caption{\label{Fig_2}
Circles in a discrete world governed by the $\ell^1$-norm. Using Euclid's original definition of a circle as the set of all points having the same distance from a given center, a ``circle'' emerging on a regular grid is made of points which are arranged in a square-like fashion. If we increase the radius $r$ of the circle, the set of points forming our ``circle'' takes the shape of a square standing on one of its corners (right).}
\end{figure}

But the surprises we can find in Manhattan don't stop there. What about geometrical figures in a space, discrete or continuous, endowed with $\ell^1$-norm? Let us consider a circle, for example. Being one of the most simple, and most-studied, geometrical figures, we follow Euclid's original definition, which sees a circle as \emph{``a plane figure contained by one line such that all the straight lines falling upon it from one point among those lying within the figure are equal to one another.''} (Euclid's \emph{Elements}, Book I, \S 19). If we adhere to this definition and consider lines as equal to one another if their $\ell^1$-norms are equal, we end up with an arrangement of points forming our circle which, by no means, resembles that of the circle we are all too familiar with (Fig.~\ref{Fig_2}). Even if we increase the size of this ``discrete circle'', and ultimately reach the continuum limit, we remain with the shape of a square rotated by $\pi/4$. Still more amusing, let us calculate pi in this ``circle''. If we define the radius as the distance $r$ between the center $O$ and each of the corners $A, B, C$ and $D$, then, according to the definition of the Manhattan distance (\ref{Eq_l1Norm}), the four points $A, B, C$ and $D$ are spread $2r$ away from each other, yielding a circumference of $8r$ for our circle. As pi is the constant which defines the ratio between circumference and diameter of a circle, we end up with $\pi = 4$, a delightful integer value which, however, is also a far shot off the  transcendental value of about 3.14159 each mathematics teacher impresses his or her young pupils with. 

Of course one could argue that the case considered above, an Euclidean plane endowed with $\ell^1$-norm, is of mere academic interest without actual link to reality, apart from the reality faced each day by taxicab drivers and their passengers in Manhattan. But we could also turn the table around and ask, in the best tradition of mathematics: what if such a discrete space with its $\ell^1$-norm is real? That this question is not too far-fetched is supported by the fact that, since many decades, vast numbers of mathematical physicists indulge in a desperate search for a quantum theory of gravity \cite{Rovelli07,Merali13}, a valid model of a discrete spacial (and temporal) makeup of reality itself. Together with healthy philosophical elaborations \cite{Hagar14} and emerging experimental results (e.g., \cite{AlbertEA08}, but see \cite{Nemiroff12}), the consequences of such a possibility are far-reaching because, eventually or ultimately, they might lead to the rejection of the ideal real number line in favour of a discrete and finite (or effinite, see \cite{Gauthier02}) mathematical underpinning of the very construct of reality. But what does this mean for our squared ``circle''?


\section{Getting round in a discrete world.}

As we saw above, when using Euclid's original circle definition, the resulting geometrical figure in our discrete world carries no resemblance to a circle in $\mathbb{R}^2$, even if we approach the asymptotic limit of continuity. However, Euclid's definition is just one of several. A mathematically perhaps more rigorous definition is given by a parametric representation, which defines a circle $S^1 \subset \mathbb{R}^2$ as the set of all points $(x,y) \in \mathbb{R}^2$ which satisfy the algebraic relation
\begin{equation}
\label{Eq_circleP}
(x-x_0)^2 + (y-y_0)^2 = r^2,
\end{equation}
where $r \in \mathbb{R}: r > 0$ is called the radius of the circle and the point $(x_0,y_0)$ its origin. 

Inspired by this representation, we will construct now a set of points which not only describes a circle in our discrete world with $\ell^1$-norm, but also, if we make the circle sufficiently big, does take on the shape we are all familiar with. For that, we need to enter the field of ``digital geometry'', a branch of mathematics which deals with the digital, hence discrete, representation of geometrical figures \cite{KletteRosenfeld04}. Driven by technological applications, this field emerged together with digital devices, such as plotters and matrix displays, in the middle of the last century. Even modern computer screens, although their regular grid  layout may not be visible to anyone anymore, rely on digital geometric drawing primitives. 

\begin{figure}[t!]
\centering
\includegraphics[width=\textwidth]{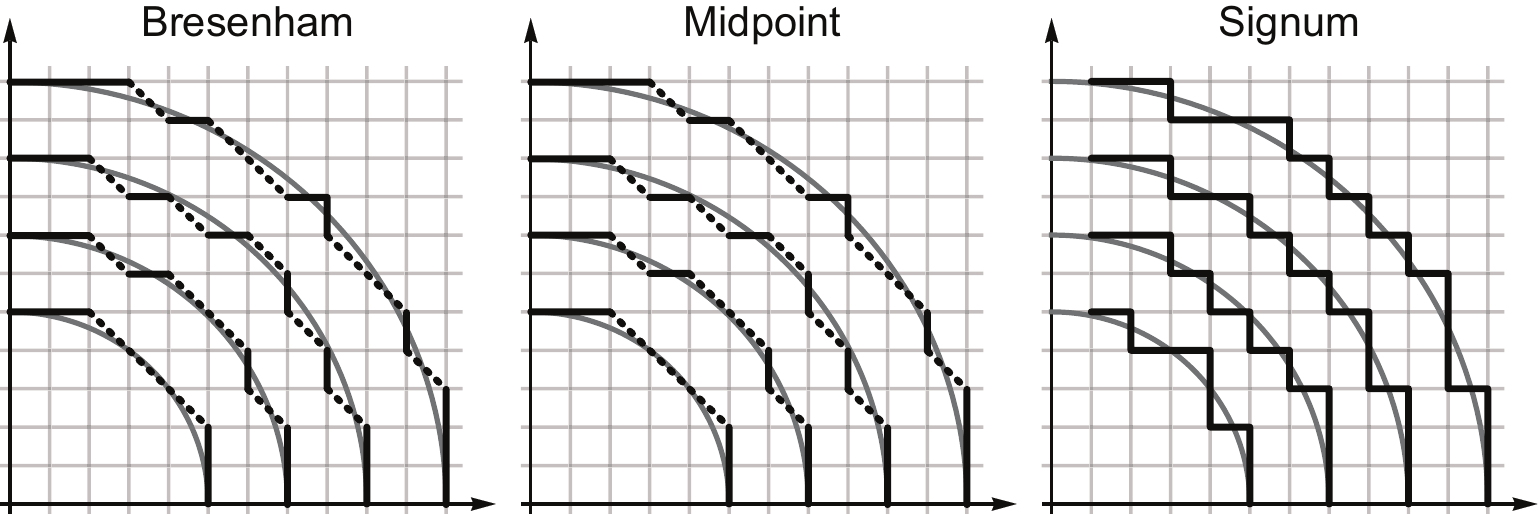}
\caption{\label{Fig_3}
Construction of circles on an integer lattice $\mathbb{Z}^2$ using different algorithms (Bresenham \cite{Bresenham77}; Midpoint \cite{FoleyEA90}; Signum: see text and \cite{Rudolph16a}; for a thorough comparative study and some historical notes, see \cite{BarreraEA15}). Shown are examples of digital circles (black) approximating Euclidean circles of radii 5, 7, 9 and 11 (gray). With the exception of the signum algorithm, most of the digital circle algorithms cited in the literature do not yield valid paths on $\mathbb{Z}^2$ (black dotted; see \cite{Rudolph16a}).}
\end{figure}

A great number of such primitives exist for curved shapes, circles in particular, with the Bresenham \cite{Bresenham77} and Midpoint \cite{FoleyEA90} algorithms being the most widely used (Fig.~\ref{Fig_3}). A common feature among most of these algorithms is that they do construct approximations of circles by utilizing vertical, horizontal and sloped lines. In other words, neighbouring points delivered by these algorithms are typically not neighbouring points on a regular grid. Moreover, due to the focus on speedy drawing primitives for technical applications, these algorithms are typically implemented computationally, and rely on decision trees and case distinctions. For a mathematical study of digital geometrical shapes, however, a rigorous algebraic implementation is required. 

\begin{figure}[t!]
\centering
\includegraphics[width=\textwidth]{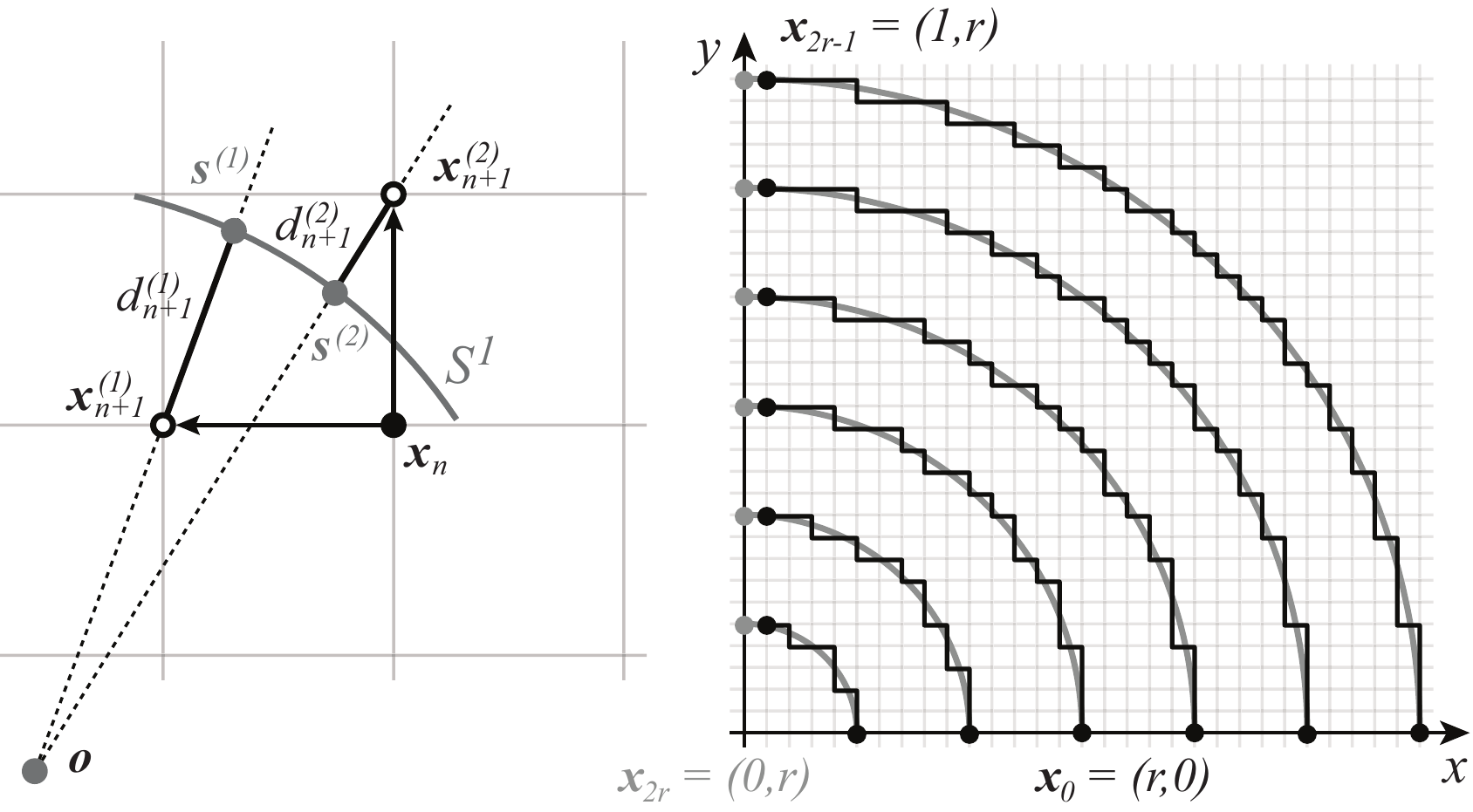}
\caption{\label{Fig_4}
Recursive algebraic construction of a circular path $\mathcal{S}^1$ on $\mathbb{Z}^2$ in the upper right quadrant, approximating $S^1 \subset \mathbb{R}^2$ (left), and examples of digital circles of various integer radii ($r=5,10,15,20,25,30$) constructed by using the signum algorithm (right).}
\end{figure}

In order to construct algebraically a circular path $\mathcal{S}^1$ on the regular integer lattice $\mathbb{Z}^2$ which approximates a circle $S^1$ of integer radius $r$ in $\mathbb{R}^2$, we will restrict for reasons of symmetry and notational simplicity to the upper right quadrant, and assume the origin of the circle $\boldsymbol{o} = (0,0)$. If we denote by $\boldsymbol{x}_n = (x_n,y_n) \in \mathcal{S}^1$ with $x_n, y_n \in \mathbb{Z}, n \in \mathbb{N}$ a point already assigned to the circular path we are looking for, then there are only two possibilities for the next, yet unassigned neighbouring point $\boldsymbol{x}_{n+1}$, namely
\begin{equation}
\label{Eq_xn1}
\boldsymbol{x}_{n+1} = (x_{n+1},y_{n+1}) =
\left\{
\begin{array}{l}
\boldsymbol{x}_{n+1}^{(1)} = (x_n-1,y_n), \\[0.2em]
\text{or, } \boldsymbol{x}_{n+1}^{(2)} = (x_n,y_n+1)
\end{array}
\right.
\end{equation}
(see Fig.~\ref{Fig_4}, left). For deciding between $\boldsymbol{x}_{n+1}^{(1)}$ and $\boldsymbol{x}_{n+1}^{(2)}$, we follow an approach similar to that used in most of the known digital circle algorithms, namely utilizing the minimization of a ``cost function''. To construct the latter, we consider the intersections $\boldsymbol{s}^{(1)}$, $\boldsymbol{s}^{(2)}$ on $S^1$ of lines through $\boldsymbol{o}$ and $\boldsymbol{x}_{n+1}^{(1)}$, $\boldsymbol{x}_{n+1}^{(2)}$, respectively. The emerging line segment $\overline{\boldsymbol{s}^{(1)} \boldsymbol{x}_{n+1}^{(1)}}$ has an Euclidean length of
\begin{equation}
\label{Eq_dn11}
d_{n+1}^{(1)} 
= \left| r - \sqrt{(x_n-1)^2 + y_n^2} \right|
= \left| r - \sqrt{\frac{1}{2} a_n^2 + \frac{1}{2} b_n^2 - a_n - b_n + 1} \right| ,
\end{equation}
and the line segment $\overline{\boldsymbol{s}^{(2)} \boldsymbol{x}_{n+1}^{(2)}}$ a length of
\begin{equation}
\label{Eq_dn12}
d_{n+1}^{(2)} 
= \left| r - \sqrt{x_n^2 + (y_n+1)^2} \right| 
= \left| r - \sqrt{\frac{1}{2} a_n^2 + \frac{1}{2} b_n^2 + a_n - b_n + 1} \right| ,
\end{equation}
where $a_n = x_n + y_n$ denotes the Manhattan distance of a point $\boldsymbol{x}_{n}$ in the upper right quadrant to the center $\boldsymbol{o}$, and $b_n = x_n - y_n$. With this, we now choose $\boldsymbol{x}_{n+1}^{(1)}$ or $\boldsymbol{x}_{n+1}^{(2)}$ depending on which of their associated line segments is shortest, i.e.  
\begin{equation}
\label{Eq_xn1a}
\boldsymbol{x}_{n+1} =
\left\{
\begin{array}{ll}
\boldsymbol{x}_{n+1}^{(1)} & \text{ if } d_{n+1}^{(1)} \leq d_{n+1}^{(2)} \\[0.2em]
\boldsymbol{x}_{n+1}^{(2)} & \text{ if } d_{n+1}^{(1)} > d_{n+1}^{(2)} .
\end{array}
\right.
\end{equation}

In order to arrive at an algebraic formulation of this minimization criterion and its associated cost function, we utilize the signum, or sign, function, defined here, in slight deviation from the commonly used notion, but without loss of generality, as 
\begin{equation}
\label{Eq_Sgn}
\sgn(x) =
\left\{
\begin{array}{ll}
-1 & \text{ if } x \leq 0 \\
1  & \text{ if } x > 0 .
\end{array}
\right.
\end{equation}
Denoting with
\begin{equation}
\label{Eq_Deltan}
\Delta_n = d_{n+1}^{(1)} - d_{n+1}^{(2)}
\end{equation}
the difference in length of both line segments, (\ref{Eq_xn1a}) yields in the case of $d_{n+1}^{(1)} \leq d_{n+1}^{(2)}$, i.e. $\sgn(\Delta_n)=-1$, for the $x$-coordinate of the next point on the circular path $x_{n+1} = x_{n+1}^{(1)} = x_n-1$, and in the case of $d_{n+1}^{(1)} > d_{n+1}^{(2)}$, i.e. $\sgn(\Delta_n)=1$, the new $x$-coordinate $x_{n+1} = x_{n+1}^{(2)} = x_n$. This can be expressed algebraically, using the signum function, in form of
\begin{equation}
x_{n+1} = \frac{1}{2} ( 1 - s_n ) ( x_n - 1 ) + \frac{1}{2} ( 1 + s_n ) x_n = x_n + \frac{1}{2} \sgn(\Delta_n) - \frac{1}{2} .
\end{equation}
Similarly, for the $y$-coordinate of the next point along the circular path, we obtain
\begin{equation}
y_{n+1} = \frac{1}{2} ( 1 - s_n ) y_n + \frac{1}{2} ( 1 + s_n ) ( y_n + 1 ) = y_n + \frac{1}{2} \sgn(\Delta_n) + \frac{1}{2} .
\end{equation}

What remains is to define the first point on $\mathcal{S}^1$, and to somewhat simplify the argument $\Delta_n$ of the signum function. As we restricted in this recursive construction of a circular path on $\mathbb{Z}^2$ to the upper right quadrant, the initial point $\boldsymbol{x}_0$ will lie on the horizontal axis, i.e. $\boldsymbol{x}_0 = (r,0)$. Furthermore, by construction, the full quarter circle will be reached for $n=2r$, with $\boldsymbol{x}_{2r} = (0,r)$ residing on the vertical axis (Fig.~\ref{Fig_4}). 

With respect to the simplification of $\Delta_n$, we observe that the signum function, as operator mapping the real number line into a Boolean set, i.e. $\sgn(x): \mathbb{R} \rightarrow \{-1,1\}$, has some interesting properties. Specifically, in the case of positive $x, y \in \mathbb{R}$, we have
\begin{equation}
\sgn(x-y) = \sgn( f(x) - f(y) )
\end{equation}
for any strict monotonically increasing function $f(x): \mathbb{R} \rightarrow \mathbb{R}$. This certainly holds true for $f(x) = x^2, x \geq 0$, and allows to rewrite $\Delta_n$ in (\ref{Eq_Deltan}) as 
\begin{equation}
\label{Eq_CostSimplified}
\Delta_n = -2\left( a_n + \frac{r}{\sqrt{2}} \left( \sqrt{(a_n-1)^2+c_n^2} - \sqrt{(a_n+1)^2+c_n^2} \right) \right),
\end{equation}
where $c_n = \sqrt{b_n^2 - 2b_n + 1} = r-n-1$. With this, we finally can formulate the following

\begin{proposition}
\label{Prop_SignumAlgorithm}
A circular path $\mathcal{S}^1 \subset \mathbb{Z}^2$ approximating a circle $S^1 \subset \mathbb{R}^2$ with radius $r \in \mathbb{N}$ and origin $\boldsymbol{o} = (0,0)$ in the upper right quadrant is a set $\{ \boldsymbol{x}_n \}$ of $2r$ points $\boldsymbol{x}_n = (x_n,y_n)$ with $x_n, y_n \in \mathbb{Z}$ obeying the algebraic recursions
\begin{equation}
\label{Eq_S1Algorithm}
\left\{
\begin{array}{l}
x_0 = r , x_{n+1} = x_n + \dfrac{1}{2} s_n - \dfrac{1}{2} \\[0.7em]
y_0 = 0 , y_{n+1} = y_n + \dfrac{1}{2} s_n + \dfrac{1}{2} ,
\end{array}
\right.
\end{equation}
where $n \in [0,2r-1], n \in \mathbb{N}$, and
\begin{equation}
\label{Eq_Cost}
s_n = -\sgn\left(a_n + \frac{r}{\sqrt{2}} \left( \sqrt{(a_n-1)^2+c_n^2} - \sqrt{(a_n+1)^2+c_n^2} \right)\right)
\end{equation}
with $a_n = x_n + y_n$ and $c_n = r-n-1$.
\end{proposition}

It is important to note that, in contrast to many other known digital circle algorithms, the above recursion constructs a circular path on which two neighbouring points are also neighbours on the underlying lattice, i.e. are separated by exactly one edge or link. Such a path we will call a ``valid'' path, and it is made up of only vertical and horizontal lines (Fig.~\ref{Fig_3}, right). Furthermore, as already mentioned above, the total number of points on the circular path increases linearly with the radius $r$ of the circle. Specifically, for a quarter circle, we need to recursively construct $2r$ points, i.e. $8r$ points for the full circle, contrasting other digital circle algorithms in which only about 70\% of the points required for a valid path are delivered. 

Finally, for $r > 4$, the cost function $s_n$ can be approximated faithfully by \cite{Rudolph16a}
\begin{equation}
\label{Eq_CostApproximated}
s_n = -\sgn\left( a_n^2 + c_n^2 + 1 - 2r^2 \right).
\end{equation}
Interestingly, both in (\ref{Eq_Cost}) and its approximation (\ref{Eq_CostApproximated}), only $a_n$, i.e. the Manhattan distance of the actual point $\boldsymbol{x}_n$ to the circle's center, appears. Thus, a valid circular path $\mathcal{S}^1$ on the integer lattice $\mathbb{Z}^2$ is a set of points obeying the parametric algebraic recursion (\ref{Eq_S1Algorithm}) in the Manhattan distance, and is thus akin to the parametric definition of the circle $S^1$ on $\mathbb{R}^2$ in (\ref{Eq_circleP}). Our taxicab driver would have no problems to show the passenger ``round'' in Manhattan, performing some easy calculations at each corner to proceed with the correct path. The most pressing question remaining for our mathematically inclined passenger is, whether it is possible to also find $\pi$ along this path, using solely the count of traversed blocks along the way, and if so, how. 


\section{The search for Pi.}

\begin{figure}[t!]
\centering
\includegraphics[width=\textwidth]{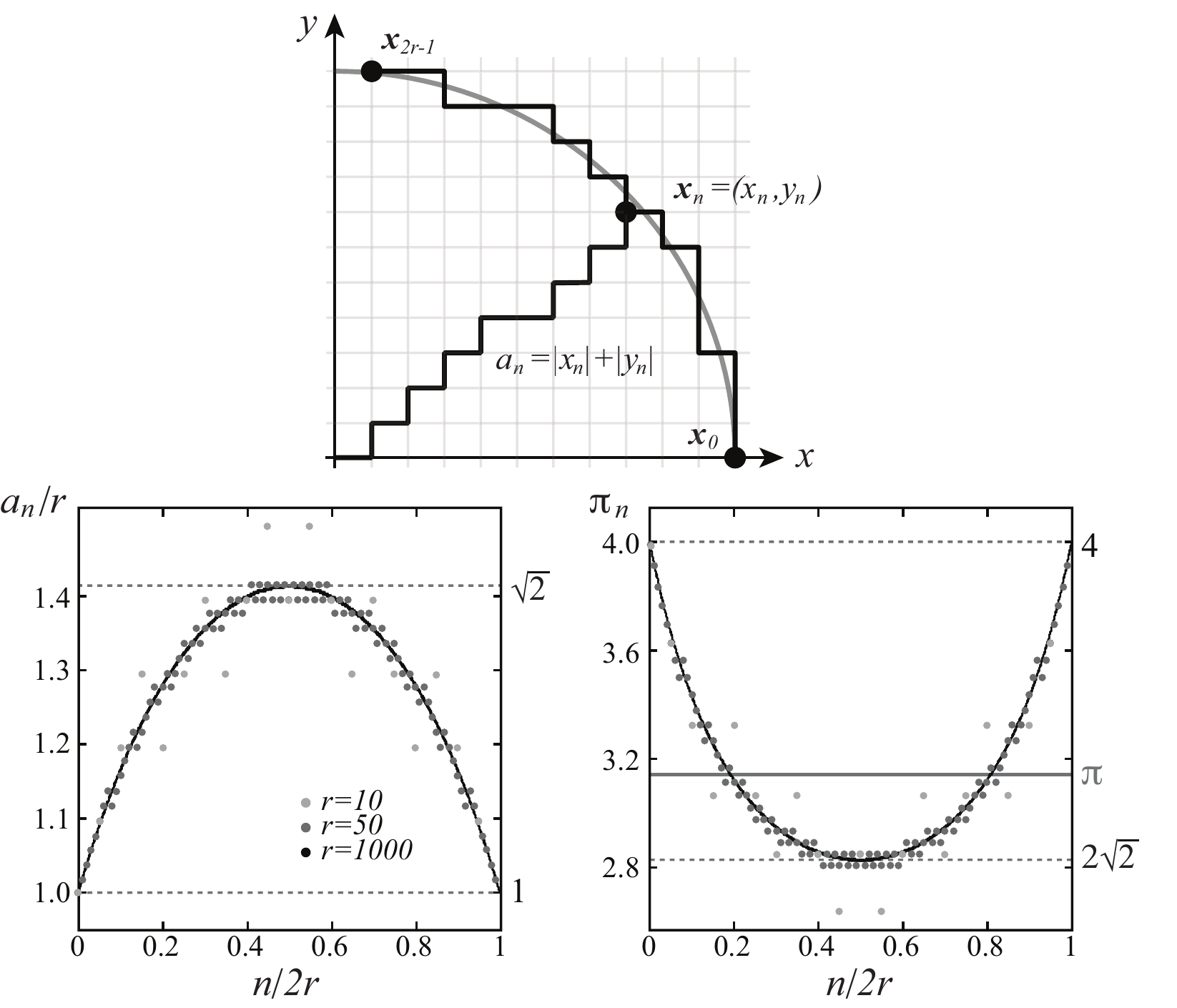}
\caption{\label{Fig_5}
The search for pi in a discrete space endowed with $\ell^1$-norm. Each point $\boldsymbol{x}_n$ along a circular path $\mathcal{S}^1 \subset \mathbb{Z}^2$ has a certain Manhattan distance $a_n$ to the circle's center which varies as we proceed along $\mathcal{S}^1$ (top). Drawing circles with larger and larger radii $r$, we find that the distance $a_n$ will progressively be bound between $r$ and $2 \sqrt{2}r$ (bottom left). Defining pi as the ratio between the total length of our path, i.e. the circumference of the circle $\mathcal{C} = 8r$, and the distance to the center, each point along $\mathcal{S}^1$ will have its own value $\pi_n$ bound, for large radii, between $2\sqrt{2}$ and 4 (bottom right).
}
\end{figure}

As we are still in Manhattan and, thus, are bound to the regular grid of streets and boulevards, we already know that, if we travel along the streets while keeping the same number of blocks between us and the virtual center of our path, we will traverse an oddly square-shaped ``circle'' (Fig.~\ref{Fig_2}) which gives rise to an integer pi of 4. On the other hand, if we follow the rules laid down in Proposition~\ref{Prop_SignumAlgorithm}, we approach a more circular path (Fig.~\ref{Fig_4}), but our distance to the virtual center will change at each corner reached. In fact, after having completed a quarter circle, thus visited $2r$ corners, we find that, for sufficiently large circles, this distance will have varied between about $r$ and $\sqrt{2}r$, as shown in Fig.~\ref{Fig_5} (bottom left). 

Being hopeful to make sense of this in our search for pi, we then calculate the ratio between the circumference of our circular path, $\mathcal{C} = 8r$, and the measured distance $a_n$ to the center at each corner. Unsurprisingly, also here we obtain values which vary along the path, and for large circle radii $r$ are progressively bound between $2\sqrt{2}$ and 4 (Fig.~\ref{Fig_5}, bottom right). Although being still somewhat puzzled by this observation, it is relieving to note that the true value of $\pi$ resides comfortably in between these two bounds. 

\begin{figure}[t!]
\centering
\includegraphics[width=\textwidth]{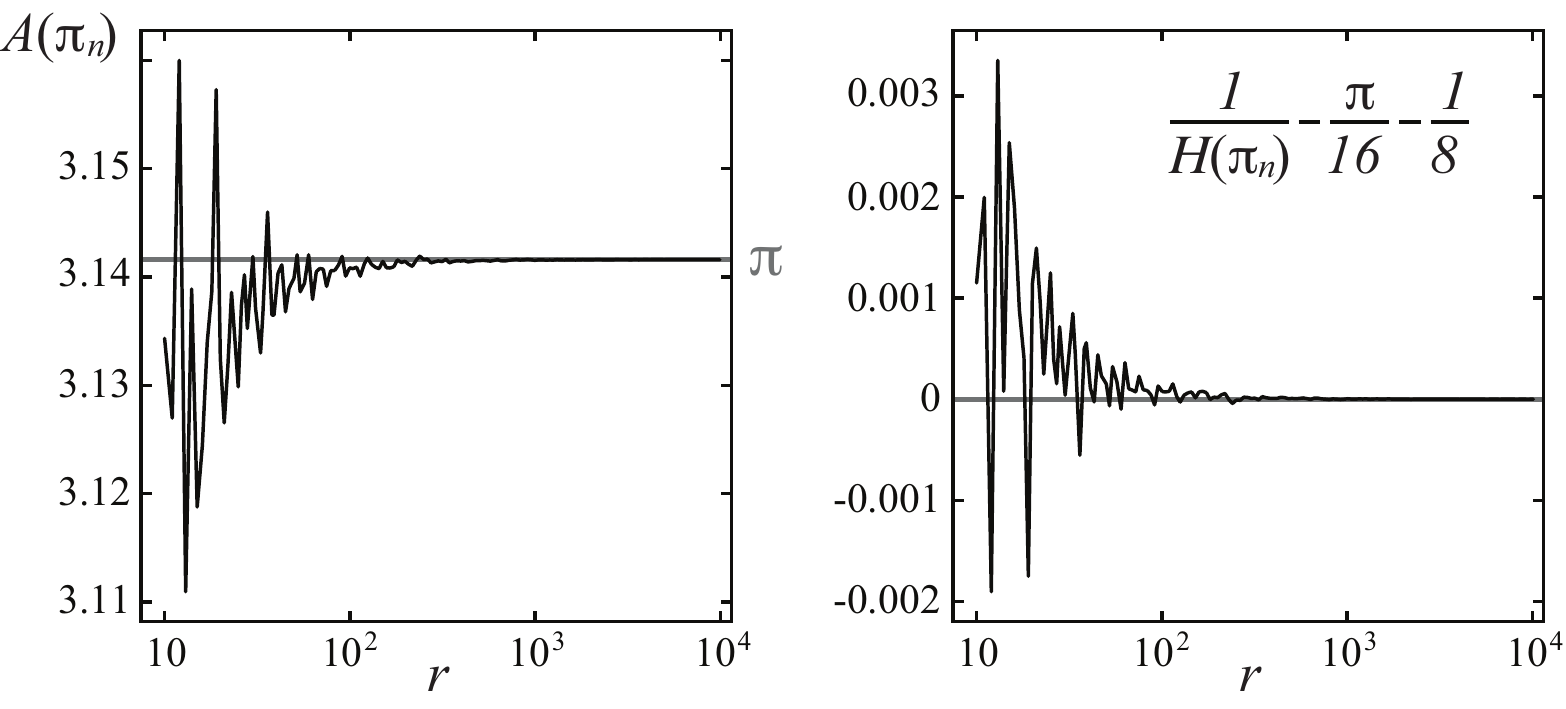}
\caption{\label{Fig_6}
Emergence of $\pi$ on the circular path $\mathcal{S}^1 \subset \mathbb{Z}^2$ with increasing radius $r$. Both the arithmetic mean $A(\pi_n)$ and harmonic mean $H(\pi_n)$ of the sequence $\pi_n$, Equation (\ref{Eq_pin}), generated recursively through (\ref{Eq_S1Algorithm}), yield $\pi$ in the asymptotic limit $r \rightarrow \infty$.}
\end{figure}

Having meticulously recorded all distances and associated pi values at each turn, we can now start to play around with the list of numbers. As already noted above, the values on this list wiggle around the familiar $\pi$, and, considering a quarter circle, come close to $\pi$ after about 19.2\% (for around $n=0.38r$) and 80.8\% ($n=1.62r$) of the way along the circular path. Following our mathematical intuition, we consider taking the arithmetic mean of the recorded numbers, and are rewarded by a somewhat baffling numerical surprise: the average value of all the recorded pi values is close to $\pi$! In fact, the larger we take our circle in Manhattan, the smaller the difference to the well-known value of $\pi$ becomes. While for a radius of $r=10$ the two handful of recorded numbers on our list average to 3.13433, traversing along a circular path of radius $r=1,000,000$ delivers $\pi$, with 3.141592652, precisely up to 9 decimals (Tab.~\ref{Tab_1} and Fig.~\ref{Fig_6} left).

\begin{table}[t!]
\centering
\caption{\label{Tab_1}
The emergence of $\pi$ in a space endowed with $\ell^1$-norm. Drawing circles with increasing radii $r$, the arithmetic mean $A(\pi_n)$ of the pi-values $\pi_n$ associated with each point $\boldsymbol{x}_n$ along the circular path $\mathcal{S}^1 \subset \mathbb{Z}^2$ approaches slowly $\pi$.}
\begin{tabular}{ | r | l | l | }
\hline
$r$ & $A(\pi_n)$ & $|A(\pi_n)-\pi|$ \\ \hline
10				& 3.134332334  & 0.007260319 \\ 
100 			& 3.140836098  & 0.000756555 \\
1,000 		& 3.141531737  & 0.000060916 \\
10,000 		& 3.141591709	& 9.45037$\cdot 10^{-7}$ \\
100,000 		& 3.141592625	& 2.87631$\cdot 10^{-8}$ \\
1,000,000	& 3.141592652	& 1.58979$\cdot 10^{-9}$ \\
\hline
\end{tabular}
\end{table}

After a few more tests and nightly trips on circular paths with larger and larger radii through the city which never sleeps, we finally feel confident enough to formulate the following

\begin{conjecture}
\label{Conj_Api}
The arithmetic mean
\begin{equation}
A(\pi_n)
= \frac{1}{2r} \sum\limits_{n=0}^{2r-1} \frac{\mathcal{C}}{2a_n(r)}
= 2 \sum\limits_{n=0}^{2r-1} \frac{1}{a_n(r)},
\end{equation}
of the finite sequence
\begin{equation}
\label{Eq_pin}
\pi_n(r) = \frac{\mathcal{C}}{2a_n(r)} = \frac{4r}{a_n(r)},
\end{equation}
where $a_n = |x_n| + |y_n|$ denotes the $\ell^1$-norm of each point $(x_n,y_n) \in \mathbb{Z}^2$ on the circular path $\mathcal{S}^1 \subset \mathbb{Z}^2$ constructed recursively by (\ref{Eq_S1Algorithm}), converges to $\pi$ in the asymptotic limit $r \rightarrow \infty$, i.e.
\begin{equation}
\lim_{r \rightarrow \infty} A(\pi_n) = \pi.
\end{equation}
\end{conjecture}

We, indeed, found $\pi$ in Manhattan! But, while being quite happy with this unexpected discovery, we cannot find rest. As mathematicians, we must set out to more rigorously prove the validity of our assertion, without the need to drive long distances and keep track of long lists of numbers at each turn.  


\section{Capturing Pi, but not quite.}

So the formal hunt for $\pi$ in Manhattan begins. Unfortunately however, we soon realize that in order to show how this upmost intriguing constant emerges out of the average on the regular grid of streets, we will be forced to leave Manhattan, and incorporate the continuous Euclidean plane on which this discrete construct lives. This does render the proof of Conjecture \ref{Conj_Api}, at least for now, less satisfying. However, in defence of the hunter, we assume that each analytical expression that we will encounter in what follows has to be understood in terms of its power expansion, so remains finite for all practical purposes. As we are interested in the asymptotic limit of a circular path with infinite radius, we then only have to show that the terms which mark the difference between the finite reality and infinite ideal will vanish as we increase the size of our discrete circle. 

The main difficulty we face in capturing $\pi$ is the occurrence of the signum function (\ref{Eq_Cost}) in the algebraic recursion (\ref{Eq_S1Algorithm}). To illustrate this point, we will first deduce the explicit form of the Manhattan distance $a_n$ of each point along the circular path. Restricting again to the upper right quadrant, we immediately find from the recursions (\ref{Eq_S1Algorithm}) that the explicit form of the coordinates along $\mathcal{S}^1$ is given by
\begin{eqnarray*}
x_n & = & r + \frac{1}{2} S_{n-1} - \frac{n-1}{2} \\
y_n & = & \frac{1}{2} S_{n-1} + \frac{n-1}{2} ,
\end{eqnarray*}
where $n \in [1, 2r-1]$ and $x_0=r, y_0=0$. Here, we introduced the partial sum
\begin{equation}
\label{Eq_Sn}
S_n = \sum\limits_{k=0}^{n} s_k 
\end{equation}
of the signum terms (\ref{Eq_Cost}). It immediately follows that the Manhattan distance $a_n = x_n + y_n$ of each point on the circular path in the upper right quadrant to the center takes the explicit form
\begin{equation}
\label{Eq_an}
a_n = r + S_{n-1}
\end{equation}
with $n \in [1, 2r-1]$ and $a_0 = r$. Thus, $a_n$ is completely determined by the sequence of partial sums of the signum terms $s_n$. Although the rather peculiar properties of $\sgn(x)$ helped us earlier to simplify the generating algorithm for our circular path, see (\ref{Eq_CostApproximated}), the lack of a useable finite analytical representation of the signum function makes it extremely difficult to further treat $s_n$ and its partial sums. 

\begin{figure}[t!]
\centering
\includegraphics[width=\textwidth]{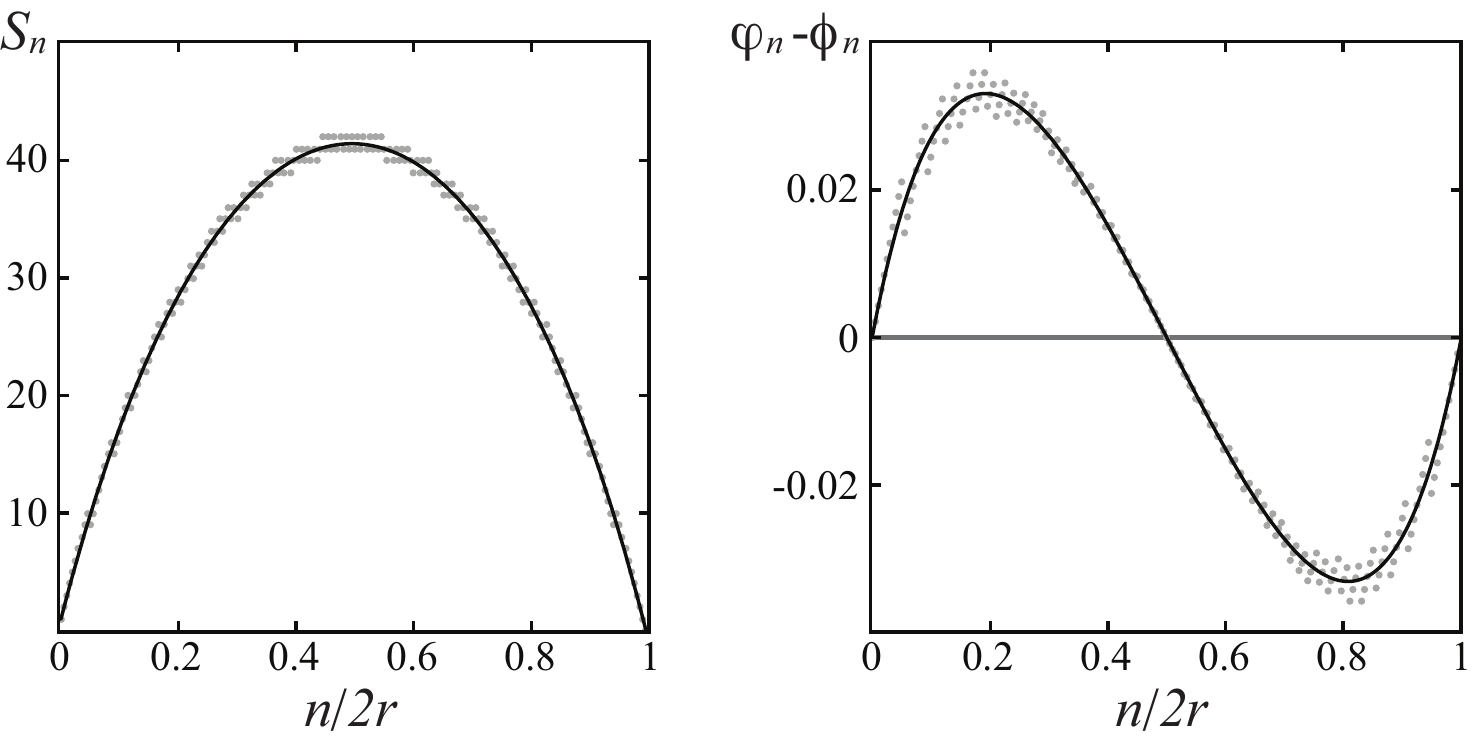}
\caption{\label{Fig_7}
The sequence of partial sums $S_n$ of the signum terms $s_n$, see (\ref{Eq_Sn}), and the absolute difference of angles $\varphi_n$, see (\ref{Eq_phin}), associated with each point along a circular path $\mathcal{S}^1 \subset \mathbb{Z}^2$ and $\phi_n = n \pi/4r$. A representative example of a circular path with radius $r=100$ (grey dots), and the analytical forms of $S_n$, Equation (\ref{Eq_SnEx}), and $\varphi_n - \phi_n$ are shown.}
\end{figure}

If we take a closer look at $S_n$ itself, we find that it takes roughly the shape of an inverted parabola, reaching its maximum for $n=r$ (Fig.~\ref{Fig_7}). Furthermore, due to its definition (\ref{Eq_Sn}), $S_n$ is subject to the recursion
\begin{equation}
\label{Eq_SnRec}
S_0 = s_0 = 1, S_{n+1} = S_n + s_{n+1}
\end{equation}
with $n \in [0, 2r-2]$, which yields $S_{2r-2} = S_{2r-1} - s_{2r-1}$. Due to symmetry of the lower and upper half of the quarter circle, the number of steps to the left ($s_n = -1$) and upwards ($s_n = 1$) must, by construction, be equal, hence $S_{2r-1} = 0$. As $s_{2r-1} = -1$, we also have $S_{2r-2} = 1$. Generalizing this argument from symmetry, other terms in the sequence $S_n$ can be treated similarly, yielding
\begin{equation*}
S_n = S_{2r-1-(n+1)}
\end{equation*}
$\forall n \in [0,r]$. Unfortunately however, these properties help little in our search for a more explicit representation of $a_n$ in (\ref{Eq_an}).

In order to proceed, we recall that we are only interested in the asymptotic limit $r \rightarrow \infty$. By construction, our discrete circular path $\mathcal{S}^1 \subset \mathbb{Z}^2$ will, with increasing precision, capture the circle $S^1 \subset \mathbb{R}^2$. Specifically, as the Euclidean distance of each point $\boldsymbol{x}_n$ along $\mathcal{S}^1$ from the continuous ideal $S^1$ remains bound, we have
\begin{equation*}
\lim_{r \rightarrow \infty} \frac{d^{(1)}_{n}}{r} = \lim_{r \rightarrow \infty} \frac{d^{(2)}_{n}}{r} = 0 .
\end{equation*}
Here we used the fact that $r \leq a_n \leq \sqrt{2} r$ for all $n \in [0,2r-1]$, and that $a_n$ reaches its maximum for $n=r$, for which $b_n = 0$. Thus, the relative Euclidean distance of each $\boldsymbol{x}_n$ to the circle's center, $\lVert \boldsymbol{x}_n \rVert_2/r$, will converge to 1 for $r \rightarrow \infty$, i.e.
\begin{equation*}
\lim_{r \rightarrow \infty} \frac{\sqrt{x_n^2+y_n^2}}{r} = 1 .
\end{equation*}
With this, we can safely introduce an angle $\varphi_n$ associated with each point $\boldsymbol{x}_n \in \mathcal{S}^1$ according to 
\begin{eqnarray}
\label{Eq_phindef}
x_n & = & r \cos(\varphi_n) \nonumber \\
y_n & = & r \sin(\varphi_n) .
\end{eqnarray}
Considering now $b_n = x_n-y_n=r-n$, we then obtain
\begin{equation*}
n = r \big( 1 + \sin(\varphi_n) - \cos(\varphi_n) \big) ,
\end{equation*}
which associates with each angle $\varphi_n$ the sequence index $n$. Noting that we are only interested in solutions covering the upper right quadrant, for which $0 \leq \varphi_n \leq \pi/2$, we obtain from this
\begin{equation}
\label{Eq_phin}
\varphi_n = \arctan \left( \frac{ nr - r^2 + r \sqrt{r^2+2nr-n^2} }{ r^2 - nr + r \sqrt{r^2+2nr-n^2} } \right),
\end{equation}
where $n \in [0,2r-1]$.

Equation (\ref{Eq_phin}) associates with each $n$ along the circular path an angle $\varphi_n$, marking the discretization of the Euclidean circle $S^1$ within our approach. However, these angles are not uniformly spread apart, but show a very distinct deviation from an equal spacing $\phi_n = n \pi / 4r$. In fact, as it was demonstrated in \cite{Rudolph16a}, it is exactly this deviation from a discretization of $S^1$ using $\phi_n$ instead of $\varphi_n$ which is responsible for obtaining $\pi$, as will be shown below. 

Now that we have established a link between the sequence index $n$ and an associated angle $\varphi_n$ for each point along the circular path $\mathcal{S}^1$, we can proceed and express the Manhattan distance $a_n$ of each of these points to the center. Recalling that, in the upper right quadrant, $a_n = x_n + y_n$, using (\ref{Eq_Sn}) and (\ref{Eq_phindef}) with (\ref{Eq_phin}) yields
\begin{equation*}
S_{n-1} 
= x_n + y_n - r 
= \sqrt{ r^2 + 2nr - n^2} - r
\end{equation*}
and, thus,
\begin{equation}
\label{Eq_SnEx}
S_n = \sqrt{ r^2 + 2(n+1)r - (n+1)^2} - r,
\end{equation}
which holds for all $n$ with $1 \leq n \leq 2r-1$. Equation (\ref{Eq_SnEx}) provides an explicit representation of the partial sum of the signum terms, which, according to construction, is valid in the asymptotic limit $r \rightarrow \infty$. With this, we can then explicitly construct, according to (\ref{Eq_pin}), the sequence of pi values $\pi_n$ associated with each point along the circular path, and calculate its arithmetic mean. We obtain
\begin{equation}
\label{Eq_Apin0}
A(\pi_n) 
= \frac{2}{a_0} + 2 \sum\limits_{n=1}^{2r-1} \frac{1}{r + S_{n-1}} 
= \frac{2}{r} + 2 \sum\limits_{n=1}^{2r-1} \frac{1}{\sqrt{ r^2 + 2nr - n^2 }}.
\end{equation}
What remains is to explicitly evaluate the sum in the last equation. 

Such an evaluation is, however, not easy due to the particular form of the denominator, and requires some juggling with number theory. First, we recall the parabolic-like shape of this term, which suggests a power expansion about $n=r$. Noting that 
\begin{equation*}
r^2+2nr-n^2 = -(n-r)^2+2r^2
\end{equation*}
and that 
\begin{equation*}
(n-r)^2 \leq r^2 < 2r^2
\end{equation*}
for all $n$ with $0 \leq n \leq r$ and $r > 0$, we have
\begin{equation*}
\frac{1}{\sqrt{ r^2 + 2nr - n^2 }} = \frac{1}{\sqrt{2} r} \frac{1}{\sqrt{1 - \frac{(n-r)^2}{2r^2}}},
\end{equation*}
where $(n-r)^2/2r^2 < 1$ for all $r,n \in \mathbb{N}$ with $0 \leq n \leq 2r-1$. We can now employ Newton's generalized binomial theorem \cite{Coolidge49, Gross07}, specifically
\begin{equation*}
\frac{1}{\sqrt{1-x}} = \sum\limits_{k=0}^{\infty} {k-\frac{1}{2} \choose k} x^k,
\end{equation*}
which holds for all $x$ with $|x|<1$. With this, the arithmetic mean (\ref{Eq_Apin0}) takes the form
\begin{eqnarray}
\label{Eq_Apin1}
A(\pi_n) 
& = & \frac{2}{r} + \frac{\sqrt{2}}{r} \sum\limits_{n=1}^{2r-1} \sum\limits_{k=0}^{\infty} {k-\frac{1}{2} \choose k} \frac{(n-r)^{2k}}{2^k r^{2k}} \nonumber \\
& = & \frac{2}{r} + \frac{\sqrt{2}}{r} \sum\limits_{k=0}^{\infty} {k-\frac{1}{2} \choose k} \frac{1}{2^k r^{2k}} \sum\limits_{n=1}^{2r-1} (n-r)^{2k},
\end{eqnarray}
where we used $a_0=r$ and reordered the two sums. 

The sum of $(n-r)^{2k}$ over $n$ in (\ref{Eq_Apin1}) can easily be expressed in terms of the generalized zeta function, or Hurwitz zeta \cite{MagnusOberhettinger66}
\begin{equation*}
\zeta(s,a) = \sum\limits_{n=0}^{\infty} \frac{1}{(n+a)^s},
\end{equation*}
yielding
\begin{equation}
\label{Eq_Apin2}
A(\pi_n) = \frac{2}{r} + \frac{\sqrt{2}}{r} \sum\limits_{k=0}^{\infty} {k-\frac{1}{2} \choose k} \frac{1}{2^k r^{2k}} \big( \zeta(-2k,1-r) - \zeta(-2k,r) \big).
\end{equation}
Exploiting
\begin{equation*}
\zeta(-n,x) = - \frac{B_{n+1}(x)}{n+1}
\end{equation*}
\cite[Theorem 12.13]{Apostol95}, which holds for $n \geq 0$ and links the Hurwitz zeta to Bernoulli polynomials
\begin{equation*}
B_{n}(x) = \sum\limits_{k=0}^n \binom{n}{k} B_{n-k} x^k,
\end{equation*}
we can further simplify $A(\pi_n)$ to
\begin{eqnarray}
\label{Eq_Apin3}
A(\pi_n) 
& = & \frac{2}{r} + \frac{\sqrt{2}}{r} \sum\limits_{k=0}^{\infty} \frac{1}{2k+1} {k-\frac{1}{2} \choose k} \frac{1}{2^k r^{2k}} \big( B_{2k+1}(r) - B_{2k+1}(1-r) \big) \nonumber \\
& = & \frac{2}{r} + \frac{\sqrt{2}}{r} \sum\limits_{k=0}^{\infty} \frac{1}{2k+1} {k-\frac{1}{2} \choose k} \frac{2}{2^k r^{2k}} B_{2k+1}(r),
\end{eqnarray}
where in the last step we utilized $B_n(1-r)=(-1)^n B_n(r)$ \cite{Lehmer88}.

What remains is to simplify the sum over $k$ in (\ref{Eq_Apin3}). To that end, we first note that the Bernoulli polynomials are defined in terms of Bernoulli numbers $B_{n}$ as
\begin{equation*}
B_n(x) = \sum\limits_{l=0}^n {n \choose l} B_{n-l} x^l. 
\end{equation*}
This leaves us with 
\begin{equation}
\label{Eq_Apin4}
A(\pi_n) = \frac{2}{r} + \sqrt{2} \sum\limits_{k=0}^{\infty} \frac{1}{2k+1} {k-\frac{1}{2} \choose k} \frac{1}{2^{k-1}}
\sum\limits_{l=0}^{2k+1} {2k+1 \choose l} B_{2k+1-l} r^{l-2k-1}.
\end{equation}
To further treat this expression, we recall that we are only interested in the asymptotic limit, i.e. $r \rightarrow \infty$. In this case, all terms $r^{l-2k-1}$ for which $l-2k-1 < 0$ in the above sum converge to 0 and can be removed. In fact, upon closer inspection, the only surviving term in the second sum of (\ref{Eq_Apin4}) is that for $l=2k+1$, leaving us with
\begin{equation}
\label{Eq_limApin0}
\lim_{r \rightarrow \infty} A(\pi_n) = \sqrt{2} \sum\limits_{k=0}^{\infty} \frac{1}{2k+1} {k-\frac{1}{2} \choose k} \frac{1}{2^{k-1}},
\end{equation}
where we used $B_0 = 1$. Finally, utilizing the binomial relation
\begin{equation*}
{k-\frac{1}{2} \choose k} = \frac{1}{2^{2k}} {2k \choose k}
\end{equation*}  
\cite[Z.45]{Gould72}, (\ref{Eq_limApin0}) can be written as
\begin{eqnarray*}
\label{Eq_limApin1}
\lim_{r \rightarrow \infty} A(\pi_n) 
& = & \sqrt{2} \sum\limits_{k=0}^{\infty} \frac{1}{2k+1} \frac{1}{2^{2k}} {2k \choose k} \frac{1}{2^{k-1}} \\
& = & 4 \sum\limits_{k=0}^{\infty} \frac{1}{2k+1} \frac{1}{2^{2k}} {2k \choose k} \left( \frac{1}{\sqrt{2}} \right)^{2k+1}.
\end{eqnarray*}
Noting that the sum in the last equation is identical to the power expansion of $\arcsin(x)$ at $x=1/\sqrt{2}$, we finally arrive at 
\begin{equation}
\lim_{r \rightarrow \infty} A(\pi_n) = 4 \arcsin\left( \frac{1}{\sqrt{2}} \right) = \pi,
\end{equation}
thus suggesting the validity of Conjecture \ref{Conj_Api}.

Or does it not? As mentioned earlier, the attempt of proving Conjecture~\ref{Conj_Api} presented above is rather unsatisfactory. Firstly, in order to deal with the signum terms occurring in the recursive generation of the circular path, in particular their partial sums (\ref{Eq_Sn}), we needed to leave the discrete layout of Manhattan's streets, and consider the underlying continuous Euclidean plane on which it is built. On this plane, the larger the radius $r$ of our circle, the closer each point along the discrete path will reside with an Euclidean distance of $r$ away from the center. For a more satisfying proof of our conjecture, we should not just stay in Manhattan itself, but preferably should explicitly deal with the signum terms themselves. 

Secondly, the proposed ``proof'' requires to take the asymptotic limit $r \rightarrow \infty$ along the way, but before reaching the end. Specifically, Equation (\ref{Eq_Apin4}) can only be treated if one discards all terms which, ultimately, vanish for $r \rightarrow \infty$, before this limit is reached. Also here, for a more satisfying proof, we should remain finite and only at the very end require the asymptotic transition, along the lines suggested in the original conjecture. 


\section{Moving to Manhattan.}

If both problems which render the proposed proof of Conjecture~\ref{Conj_Api} rather unsatisfactory could be solved, we indeed would have gained much more. Similar to the recursive algebraic construction of a circular path on the integer grid presented in Proposition~\ref{Prop_SignumAlgorithm}, an algebraic recursion can be formulated which constructs a discrete path residing fully inside a given circle $S^1 \subset \mathbb{R}^2$ of radius $r$. Moreover, based on this recursion, the area enclosed by this discrete path can be obtained in a mathematically rigorous fashion, thus providing an alternative, purely recursive formulation of the well-known yet unproven Gauss Circle Problem (e.g., see \cite{MW_GCP, Huxley96}). Whether the peculiar properties of the signum function can be exploited to allow viewing this problem in a different light, or will prove to be yet another dead end in tackling it, remains to be explored. 

Linked to the area enclosed by the circular path $\mathcal{S}^1 \subset \mathbb{Z}^2$ is another interesting surprise we discover  in Manhattan. As we saw in Conjecture~\ref{Conj_Api}, the sequence of $\pi_n$ associated with each point along $\mathcal{S}^1$ is given by the sequence of reciprocals of the Manhattan distance of each of these points to the center, i.e. $1/a_n$, yielding $\pi$ when its arithmetic mean is taken in the asymptotic limit. But what about the sequence $a_n$ itself? From (\ref{Eq_pin}), we immediately have
\begin{equation}
\label{Eq_pin1}
\frac{1}{\pi_n(r)} = \frac{2a_n(r)}{\mathcal{C}} = \frac{a_n(r)}{4r}.
\end{equation}
Somewhat surprisingly, taking the arithmetic mean of the sequence $1/\pi_n$, does yield, in the asymptotic limit, $\pi$ as well, as shown in \cite{Rudolph16a}. Specifically, we can formulate

\begin{proposition}
\label{Prop_Hpi}
The arithmetic mean
\begin{equation}
\label{Eq_Arpi}
A\left(\frac{1}{\pi_n}\right) 
= \frac{1}{2r} \sum\limits_{n=0}^{2r-1} \frac{a_n(r)}{4r}
= \frac{1}{8r^2} \sum\limits_{n=0}^{2r-1} a_n(r)
\end{equation}
of the finite sequence of reciprocal $\pi_n$ values associated with each point $(x_n,y_n) \in \mathbb{Z}^2$ along a circular path $\mathcal{S}^1 \subset \mathbb{Z}^2$ constructed recursively through (\ref{Eq_S1Algorithm}), obeys in the asymptotic limit $r \rightarrow \infty$ the identity
\begin{equation}
\label{Eq_Arpi1}
\lim_{r \rightarrow \infty} A\left(\frac{1}{\pi_n}\right) = \frac{\pi}{16} + \frac{1}{8}.
\end{equation}
\end{proposition}

Yet again we found $\pi$ in Manhattan, this time by considering the average of the sequence $1/\pi_n$. Noting that the harmonic mean $H$ is the reciprocal dual of the arithmetic mean $A$ of the same sequence, we immediately have 
\begin{equation}
\label{Eq_Hpi1}
\lim_{r \rightarrow \infty} H(\pi_n) = \left( \lim_{r \rightarrow \infty} A\left(\frac{1}{\pi_n}\right) \right)^{-1} = \frac{16}{\pi+2}.
\end{equation}
Although the fundamental inequality linking the arithmetic and harmonic means of any given sequence is not violated, as
\begin{equation*}
\lim_{r \rightarrow \infty} A(\pi_n) = \pi > \frac{16}{\pi + 2} = \lim_{r \rightarrow \infty} H(\pi_n),
\end{equation*}
the results suggest that, for the specific sequence $\pi_n$ constructed along a circular path in $\mathbb{Z}^2$, an identity exist which links both means exactly. Specifically, we have
\begin{equation}
\label{Eq_AH}
\lim_{r \rightarrow \infty} \big( A(\pi_n) H(\pi_n) + 2 H(\pi_n) - 16 \big) = 0.
\end{equation}
Whether such identities linking arithmetic and harmonic means do exist for other number sequences, or for which type of number sequences the identity (\ref{Eq_AH}) holds, might be interesting questions to be explored.

The results presented here hint at some deeper number-theoretical properties of circular paths on the integer lattice, a discrete space governed by the $\ell^1$-norm. When constructed correctly, such circular paths not just morph smoothly into their continuous ideal when made larger, but also yield, somewhat surprisingly, the defining constant of the Euclidean circle, $\pi$, without surrendering the defining peculiarities of a space endowed with a linear metric. Rushing now in a taxicab across Manhattan's busy streets, we will no longer just listen to the ponderings of the cab's driver, but also wonder what other surprises we might find in the city that never sleeps.  


\section*{Acknowledgments}

Research supported in part by CNRS. The author wishes to thank JAG Willow, S Hower and CO Caine for inspiration, valuable discussions and comments.


\end{document}